\documentclass[12pt]{article}
\usepackage{amsmath,amssymb}
\newcommand{\arab}{\renewcommand{\labelenumi}{(\arabic{enumi})}}
\newcommand{\eps}{\varepsilon}
\newcommand{\isum}{\mathop{\sum}_{i=1}^n}
\newcommand{\hsum}{\mathop{\sum}_{\tau \in H}}
\newcommand{\iprod}{\mathop{\prod}_{i=1}^n}
\newcommand{\iopl}{\mathop{\oplus}_{i\in I}}
\newcommand{\hopl}{\mathop{\oplus}_{\tau \in H}}

\newcommand{\xvw}{X(W;\nu)}
\newcommand{\xbvw}{X(W;\bar{\nu})}
\newcommand{\xbvvw}{{X}(W;\bar{\nu},\nu)}

\newcommand{\xsvw}{X(W;\nu)^{st}}

\newcommand{\xxvw}{{\frak X}(W;\nu)}
\newcommand{\xxbvw}{{\frak X}(W;\bar{\nu})}
\newcommand{\xxbvvw}{{\frak X}(W;\bar{\nu},\nu)}

\newcommand{\vin}{V_{\text{in}(\tau)}}
\newcommand{\vout}{V_{\text{out}(\tau)}}
\newcommand{\lvw}{\Lambda(W;\nu)}

\newcommand{\eit}{\tilde{e_i}}
\newcommand{\fit}{\tilde{f_i}}
\newcommand{\bb}{\bar{B}}
\newcommand{\bt}{\bar{t}}
\newcommand{\bs}{\bar{s}}

\newcommand{\dimc}{\dim_{\mbox{\tiny{$\Bbb C$}}}}
\newcommand{\Homc}{\text{Hom}_{\mbox{\tiny{$\Bbb C$}}}}
\newcommand{\inn}{\text{in}}
\newcommand{\out}{\text{out}}
\newcommand{\Imm}{\text{Im}}
\newcommand{\End}{\text{End}}
\newcommand{\Ker}{\text{Ker}}
\newcommand{\Grass}{\text{Grass}}
\newcommand{\wt}{\text{wt}}

 \newtheorem{thm}{Theorem}[subsection]
 \newtheorem{prop}[thm]{Proposition}
 \newtheorem{lemma}[thm]{Lemma}
 \newtheorem{cor}[thm]{Corollary}
 \newtheorem{defn}[thm]{Definition}

 \newtheorem{ack}{Acknowledgment}

\begin{document}
 \title{Crystal bases and quiver varieties\\
{\large (Geometric construction of crystal base II)}}
 \author{Yoshihisa Saito}
 \date{ }
 \maketitle
 \begin{abstract}
  We give a crystal structure on the set of all irreducible components
  of Lagrangian subvarieties of quiver varieties. One can show that, as a 
  crystal, it is isomorphic to
  the crystal base of an irreducible highest weight representation of 
  a quantized universal enveloping algebra.
 \end{abstract}
%
%
 \renewcommand{\thefootnote}{\fnsymbol{footnote}}
 \footnote[0]{2000 {\it Mathematics Subject Classification}. 
  Primary 17B35; Secondary 14D30,16G20.}
\section{Introduction}
  Let $U_q({\frak g})$ be the quantized universal enveloping algebra
  of a Kac-Moody Lie algebra $\frak g$ with a symmetric Cartan datum
  and $U_q^-({\frak g})$ the minus part of it. In [8],[10], Lusztig
  gave a geometric realization of $U_q^-({\frak g})$ in terms of quivers
  and defined the canonical basis of it. 
  On the other hand, motivated by the study of solvable lattice models, 
  Kashiwara defined the crystal base and the global base of $U_q^-({\frak g})$
  in an algebraic way [3]. Grojnowski and Lusztig proved the global
  base is coincide with the canonical basis due to Lusztig [1]. 
  Afterwards Kashiwara and the author constructed the crystal base of   
  $U_q^-({\frak g})$ in a geometric way [7]. More precisely they defined a 
  crystal
  structure on the set of irreducible components of a Lagrangian variety
  and proved that it is isomorphic to the crystal base of $U_q^-({\frak g})$.  

  In the crystal base theory, there is another important example $B(\lambda)$
  which is the crystal base of the irreducible highest weight module of
  $U_q({\frak g})$ with a highest weight $\lambda$. The purpose of this paper
  is to construct $B(\lambda)$ in a geometrical way. Nakajima defined 
  a new family of hyper-K\"{a}hler manifolds, called quiver varieties [15],
  [16].
  He also constructed a representation of $\frak g$
  on the middle homology groups of quiver varieties. Moreover he proved
  it is isomorphic to an irreducible highest weight representation
  by using Kac's characterization of irreducible modules [2].
  In this paper we consider
  a Lagrangian subvariety of a quiver variety, following Nakajima, and
  define a crystal structure on the set of irreducible components of it.
  Instead of Kac's characterization, we give a crystal-theoretical 
  characterization of $B(\lambda)$. By the aid of it, we show the crystal
  associated to the set of irreducible components is isomorphic to 
  $B(\lambda)$.  

  In Kashiwara's algebraic 
  construction, the global base of an irreducible highest weight module of
  $U_q({\frak g})$ was also defined [3]. We remark that it has not 
  constructed yet in a geometrical way, since a representation of
  $U_q({\frak g})$ has not.   
  
 \begin{ack}
  {\em The author would like to thank Masaki Kashiwara and Toshiyuki Tanisaki
  for stimulating discussions}. 
 \end{ack}
\section{Preliminaries}
 \subsection{ }
  Let $\frak g$ be a symmetric Kac-Moody Lie algebra and $\frak t$ its
  Cartan subalgebra. Let $\{\alpha_i|i\in I\}\subset {\frak t}^*$ and
  $\{h_i|i\in I\}\subset {\frak t}$ be the set of simple roots and simple
  coroots, respectively. We normalize the non-degenerate symmetric invariant
  bilinear form ( , ) on ${\frak t}^*$ so that $(\alpha_i,\alpha_i)\in 
  {\Bbb Z}_{>0}.$ Let $P$ be the weight lattice and $P^*$ its dual lattice.
  Then $\alpha_i\in P$ and $h_i\in P^*$.

  We define $U_q({\frak g})$ as the ${\Bbb Q}(q)$-algebra generated by
  $e_i,f_i(i\in I)$ and $q^h(h\in P^*)$ with the following defining relations:
  \begin{description}
  \item[(1)] $q^h=1$ for $h=0$ and $q^{h+h'}=q^hq^{h'}$,
  \item[(2)] $q^he_iq^{-h}=q^{\langle h,\alpha_i\rangle}e_i$ and 
              $q^hf_iq^{-h}=q^{-\langle h,\alpha_i\rangle}f_i$,
  \item[(3)] $[e_i,f_j]=\delta_{i,j}(t_i-t_i^{-1})/(q_i-q_i^{-1})$
                  where $q_i=q^{(\alpha_i,\alpha_i)}$ and $t_i=q^{(\alpha_i,
                  \alpha_i)h_i}$,
  \item[(4)] $\sum(-1)^ne_i^{(n)}e_je_i^{(b-n)}=\sum(-1)^nf_i^{(n)}f_j
                  f_i^{(b-n)}=0$\\
                  where $i\ne j$ and $b=\langle h_i,\alpha_j\rangle$.
  \end{description}
 Here we used the notations $[n]_i=(q_i^n-q_i^{-n})/(q_i-q_i^{-1})$, $[n]_i!
 =\prod_{k=1}^n [k]_i$, $e_i^{(n)}=e_i^n/[n]_i!$ and $f_i^{(n)}=f_i^n/[n]_i!$.
 We understand $e_i^{(n)}=f_i^{(n)}=0$ for $n<0$. We set
 $Q=\isum {\Bbb Z}\alpha_i,\text Q_+=\isum {\Bbb Z}_{\geq 0}\alpha_i,
   \text{ and } Q_-=-Q_+$. Let $P_+$ be the set of dominant integral weights.\\
 We denote by $U_q^-({\frak g})$ the ${\Bbb Q}(q)$-subalgebra of 
 $U_q({\frak g})$ generated by $f_i(i\in I)$.

 For a fixed $i\in I$, let $U_q({\frak g}_i)$ be the ${\Bbb Q}(q)$-subalgebra 
 of $U_q({\frak g})$ generated by $e_i$, $f_i$ and $q^{\pm h_i}$.  
 We say that $U_q({\frak g})$-module $M$ is integrable if $M$ has the
 weight space decomposition $M=\mathop{\oplus}_{\nu\in P}M_{\nu}$ and 
 $M$ is $U_q({\frak g}_i)$-locally finite for any $i\in I$.  
\subsection{crystals}
 The theory of crystal base developed in [3] provides very powerful method 
 in the representation theory of $ U_q({\frak g})$. Motivated by the properties
 of the crystal base, Kashiwara defined {\it crystals} in the conbinatrial way 
 [4].
 First we recall the definition of crystals. See [3],[4],[5],[6] and [7]
 for details.

 \begin{defn}
  {\em { A crystal B is a set with}
  $$ \mbox{maps } \wt:B\to P,\mbox{ } \varepsilon_i:B\to {\bf Z}\sqcup
  \{-\infty \}\mbox{ and } \varphi_i:B\to {\bf Z}\sqcup \{-\infty \},
  \leqno{(2.2.1)}$$
  $$\tilde{e}_i :B\to B\sqcup \{0\},\mbox{ } \tilde{f}_i :B\to B\sqcup \{0\}. 
  \leqno{(2.2.2)}$$
  { They are subject to the following axioms}:
  \begin{description}
  \item[(C 1)] \(\varphi_i (b)=\varepsilon_i (b)+\langle h_i,\wt(b)\rangle.\)
  \item[(C 2)] If \(b\in B\) and \(\tilde{e}_i b\in B\) then,\\
  \(\wt(\tilde{e}_i b)=\wt(b)+\alpha_i\), \(\varepsilon_i(\tilde{e_i} b)=
  \varepsilon_i (b)-1\) and \(\varphi_i (\tilde{e_i} b)=\varphi_i(b)+1.\)
  \item[(C 2')] If \(b\in B\) and \(\tilde{f_i} b\in B\), then\\
   \(\wt(\tilde{f_i} b)=\wt(b)-\alpha_i\), \(\varepsilon_i (\tilde{f_i} b)=
   \varepsilon_i(b)+1\) and \(\varphi_i(\tilde{f_i} b)=\varphi_i(b)-1\).
  \item[(C 3)] For \(b,b'\in B\) and \(i\in I\), \(b'=\tilde{e_i} b\) if and 
   only if
  \(b=\tilde{f_i} b'\).
  \item[(C 4)] For \(b\in B\), if \(\varphi_i(b)=-\infty \), then 
  \(\tilde{e_i} b
  =\tilde{f_i} b=0.\)
  \end{description}}
  \end{defn}
  For two crystals \(B_1\) and \(B_2\), a morphism \(\psi\) from \(B_1\) to
  \(B_2\) is a map \(B_1\to B_2\sqcup \{0\}\) that satisfies the following
  conditions:
  \begin{description}
  \item[(2.2.3)] {If \(b\in B_1\) and \(\psi (b)\in B_2\), then 
                 \(\wt(\psi(b))=\wt(b)\), \(\varepsilon_i (\psi (b))=
                 \varepsilon_i (b)\), and \(\varphi_i(\psi (b))=\varphi_i 
                 (b)\),}
  \item[(2.2.4)] { For \(b\in B_1\), we have \(\psi (\tilde{e_i} b)=
                 \tilde{e_i} \psi (b)\) provided \(\psi (b)\) and
                 \(\psi (\tilde{e_i} b)\in B_2\)},
  \item[(2.2.5)] { For \(b\in B_1\), we have \(\psi (\tilde{f_i} b)=
                 \tilde{f_i} \psi (b)\) provided \(\psi (b)\) and
                 \(\psi (\tilde{f_i} b)\in B_2\)}.
  \end{description}
  A morphism \(\psi :B_1\to B_2\) is called {\it strict}, if it commutes with
  all \(\tilde{e_i}\) and \(\tilde{f_i}\).\\
  
  For two crystals \(B_1\) and \(B_2\), we define its tensor product
  \(B_1\otimes B_2\) as follows:
 \[B_1\otimes B_2=\left\{ b_1\otimes b_2~;~b_1\in B_1\mbox{ } and\mbox{ } b_2
  \in B_2\right\},\]
  \[\varepsilon_i (b_1\otimes b_2)=\max\left\{\varepsilon_i (b_1),\mbox{ }
  \varepsilon_i (b_2)-\wt_i(b_1)\right\},\]
  \[\varphi_i (b_1\otimes b_2)=\max\left\{\varphi_i (b_1)+\wt_i(b_2),\mbox{ }
  \varphi_i (b_2)\right\},\]
  \[\wt(b_1\otimes b_2)=\wt(b_1)+\wt(b_2),\]
  \[\tilde{e_i} (b_1\otimes b_2) = \left\{
  \begin{array}{rl}
   \tilde{e_i} b_1\otimes b_2 &\quad\mbox{if $\varphi_i (b_1)\geq \varepsilon_i
  (b_2)$},\\
   b_1\otimes \tilde{e_i} b_2 &\quad\mbox{if $\varphi_i (b_1)<\varepsilon_i
  (b_2)$},
  \end{array}\right. \]
  \[\tilde{f_i} (b_1\otimes b_2) = \left\{
  \begin{array}{rl}
   \tilde{f_i} b_1\otimes b_2 &\quad\mbox{if $\varphi_i (b_1)>\varepsilon_i
  (b_2)$},\\
   b_1\otimes \tilde{f_i} b_2 &\quad\mbox{if $\varphi_i (b_1)\leq \varepsilon_i
  (b_2)$}.
  \end{array}\right. \]
  Here \(\wt_i(b)=\langle h_i,\wt(b)\rangle\).\\
 \\
 {\bf Example 2.2.2}\quad
  For \(\lambda\in P_+\), \(B(\lambda )\) is the crystal associated with the
  crystal base of the simple highest weight module with highest weight 
  \(\lambda \). For $b\in B(\lambda)$ we set $\varepsilon_i (b)
  =\max\{k\geq 0~|~\tilde{e_i}^k b\ne 0\}$, $\varphi_i (b)
  =\max\{k\geq 0~|~\tilde{f_i}^k b\ne 0\}$ and $\wt(b)$ is the weight of $b$.
 We denote by $b(\lambda)$ the unique element with weight $\lambda$.\\
 \\
 {\bf Example 2.2.3}\quad
  \(B(\infty )\) is the crystal associated with the crystal base of 
  $U_q^-({\frak g})$. For $b\in B(\infty)$ we set $\varepsilon_i (b)=
  \max\{k\geq 0~|~\tilde{e_i}^k b\ne 0\}$ and $\varphi_i (b)=\varepsilon_i (b)
  +\langle h_i,\wt(b)\rangle$.\\
  \\
  {\bf Example 2.2.4}\quad
  Let $\lambda\in P_+$ be a dominant integral weight. Consider the set 
  $T_{\lambda}=\{t_{\lambda}\}$ with one element. Define $\wt(t_{\lambda})
  =\lambda$, $\eps_i(t_{\lambda})=\varphi_i(t_{\lambda})=-\infty$ and
  $\tilde{e_i}(t_{\lambda})=\tilde{f_i}(t_{\lambda})=0$ for all $i\in I$.
  Then $T_{\lambda}$ is a crystal.
\subsection{ }
 In this subsection we shall give a crystal-theoretical characterization of 
 $B(\lambda)$. We will use this result in $\S$4.  
 Consider the tensor product of crystals $B(\infty)\otimes T_{\lambda}$.
 Let $\pi_{\lambda}:
 B(\infty)\otimes T_{\lambda}\to B(\lambda)$ be a strict morphism given by
 $\tilde{f_{i_1}}
 \cdots\tilde{f_{i_l}}b_0\otimes t_{\lambda}\mapsto\tilde{f_{i_1}}
 \cdots\tilde{f_{i_l}}b(\lambda)$. The morphism $\pi_{\lambda}$ has following
 properties:
$$\begin{array}{ll}
\mbox{The set $\{b\in B(\infty)\otimes T_{\lambda}~|
 ~\pi_{\lambda}(b)\ne 0\}$ is isomorphic to $B(\lambda)$}\\ 
 \mbox{through $\pi_{\lambda}$}.
\end{array}\leqno{(2.3.1)}$$   
 \begin{prop}
  Let $B$ be a crystal and $b_{\lambda}$ an element of $B$ with weight
  $\lambda\in P_+$. Assume the following conditions.
  \begin{enumerate}
   \arab
   \item $b_{\lambda}$ is the unique element of $B$ with weight 
         $\lambda$. 
   \item There is a strict morphism $\Phi:B(\infty)\otimes T_{\lambda}\to B$
         such that $\Phi(b_0\otimes t_{\lambda})=b_{\lambda}$ and 
         $\text{Im}\Phi=B\sqcup \{0\}$. Here $b_0$ is the unique element
         of $B(\infty)$ with weight zero.
   \item Consider the set $\{b\in B(\infty)\otimes T_{\lambda}|\Phi(b)
         \ne 0\}$. Then it is isomorphic to $B$ through $\Phi$ as a set.
   \item For any $b\in B$ and $i\in I$, $\eps_i(b)=\mbox{max}\{k\geq 0|
         \eit^k(b)\ne 0\}$ and $\varphi_i(b)=\mbox{max}\{k\geq 0|
         \fit^k(b)\ne 0\}$.
      \end{enumerate}
  Then $B$ is isomorphic to $B(\lambda)$.
 \end{prop}
 {\it Proof.} By (1) and (2) any element of $B$ has the form 
 $\tilde{f_{i_1}}\cdots \tilde{f_{i_l}}
 b_{\lambda}$ with $i_1,\cdots,i_l\in I$. The following lemma is a key of 
 the proof of Proposition 2.3.1.
 \begin{lemma}
  For a given $(i_1,\cdots,i_k)\in I^k$, the following two statements 
  are equivalent.
  $$(a)\quad\tilde{f_{i_1}}\cdots \tilde{f_{i_l}}b({\lambda})=0.\qquad
  (b)\quad \tilde{f_{i_1}}\cdots \tilde{f_{i_l}}b_{\lambda}=0.$$
 \end{lemma} 
 {\it Proof of Lemma 2.3.2.} We shall only show that $(a)$ implies $(b)$. 
 For $\nu\in
 Q_-$ let $B(\lambda)_{\lambda+\nu}$ ({\it resp.} $B_{\lambda+\nu}$) be
 the set of all elements of $B(\lambda)$ ({\it resp.} $B$) with weight 
 $\lambda+\nu$. We shall prove the statement by the induction 
 on $|\mbox{ht}(\nu)|$. Here $\mbox{ht}(\nu)=\sum_{i\in I}n_i\in {\Bbb Z}$ 
 for $\nu=\sum_{i\in I}n_i\alpha_i\in Q$. If $|\mbox{ht}(\nu)|=0$, the 
 statement is 
 clear. Assume the statement holds for $|\mbox{ht}(\nu)|\leq l-1$. Take
 $\tilde{f_{i_2}}\cdots \tilde{f_{i_l}}b({\lambda})\in 
 B(\lambda)_{\lambda+\nu}$ with $|\mbox{ht}(\nu)|=l-1$ such that
 $\tilde{f_{i_1}}(\tilde{f_{i_2}}\cdots \tilde{f_{i_l}}b({\lambda}))=0.$  
 By the definition of $B(\lambda)$ we have $\varphi_{i_1}(\tilde{f_{i_2}}
 \cdots \tilde{f_{i_l}}b({\lambda}))=0$. Since $\pi_{\lambda}$ is 
 strict, we have 
 $ \tilde{f_{i_2}}\cdots \tilde{f_{i_l}}(b_0\otimes t_{\lambda})\ne 0$
 and $\varphi_{i_1}(\tilde{f_{i_2}}\cdots \tilde{f_{i_l}}(b_0\otimes 
 t_{\lambda}))=0$. By the conditions (2), (3) in Proposition 2.3.1 and the 
 induction hypothesis we have 
 $$
 \Phi(\tilde{f_{i_2}}\cdots \tilde{f_{i_l}}(b_0\otimes t_{\lambda}))=
 \tilde{f_{i_2}}\cdots \tilde{f_{i_l}}b_{\lambda}
 \ne 0
 $$
 and
 \begin{align*}
 \varphi_{i_1}(\tilde{f_{i_2}}\cdots \tilde{f_{i_l}}(b_0\otimes t_{\lambda})) 
 &= \varphi_{i_1}(\Phi(\tilde{f_{i_2}}\cdots \tilde{f_{i_l}}(b_0\otimes 
 t_{\lambda})))\\
 &= \varphi_{i_1}(\tilde{f_{i_2}}\cdots \tilde{f_{i_l}}b_{\lambda})\\
 &=0.
 \end{align*} 
 By the condition (4) in Proposition 2.3.1 this means $\tilde{f_{i_1}}
 \tilde{f_{i_2}}\cdots \tilde{f_{i_l}}b_{\lambda}=0$. Therefore the statement
 holds for $|\mbox{ht}(\nu)|=l$. \hfill$\square$\\
 
 Let us return to the proof of Proposition 2.3.1. From (3), (2.3.1) and 
 Lemma 2.3.2, we can define a bijection $\psi:B(\lambda)\to B$ by
 $\tilde{f_{i_1}}\cdots \tilde{f_{i_l}}b({\lambda})\mapsto
  \tilde{f_{i_1}}\cdots \tilde{f_{i_l}}b_{\lambda}$. By the construction
 it is easy to see that $\psi$ is an isomorphism of crystals.
 \hfill$\square$
 \section{Quiver varieties}
  \subsection{ }
  Let $A=(a_{ij})$ be the Cartan matrix of $\frak g$. We shall define an 
  oriented graph $(I,H)$ associated with $A$ as follows.
 
  Let $I$ be the set of vertices and $H$ the set of arrows of our 
 graph. For $i,j
  \in I$ $(i\ne j)$, there are $|a_{ij}|$ arrows from $i$ to $j$.  
  Let $\out(\tau)$ ({\it resp.} $\inn(\tau)$) be the outgoing
  ({\it resp.} incoming) vertex of $\tau\in H$. For $\tau \in H$, we denote
  by $\bar{\tau}$ the same edge as $\tau$ with the reverse orientation.
  The map $\bar{ }$ is a fixed free involution of $H$.     
  An orientation of our graph is a choice of a subset $\Omega \subset H$ such 
  that $ \Omega \cup \bar{\Omega} =H,$ $ \Omega \cap \bar{\Omega} =\phi .$    
  We call our oriented graph $(I,H)$ a quiver.
 \subsection{ } 
  Let $\cal V$ be the family of $I$-graded vector spaces $V=\iopl V_i$.
  For $V\in {\cal V}$, the dimension of
  $V$ is defined to be the vector $\dim V=(\dimc V_i)_{i\in I}\in 
  {\Bbb Z}_{\geq 0}^I$. 
  For $\nu\in {\Bbb Z}_{\geq 0}^I$ let ${\cal V}_{\nu}$ be the family
  of $I$-graded vector spaces $V$ with $\dim V=\nu$.
  
  Consider another $I$-graded complex vector space $W$ with $\dim W=\lambda$
  and fix it through out this paper. From now on we 
  regard the dimension vectors $\nu$ and $\lambda$ as elements of $P$
  by the following way;
  $$\nu\mapsto -\isum \dimc V_i\alpha_i,\qquad \lambda\mapsto
  \isum \dimc W_i\Lambda_i.$$ 
  Here $\Lambda_i$ is the fundamental weight of $\frak g$. 
 
   We define a complex vector space $\xvw$ by
  $$ \xvw =(\hopl \Homc (\vout ,\vin))\oplus (\iopl \Homc (V_i,W_i))
     \oplus (\iopl \Homc (W_i,V_i)).$$
 
  For an element of $\xvw$ we denote its components 
  by $( B_{\tau},t_i,s_i)$. We write $B,t,s$ for the 
  collection $(B_{\tau})_{\tau\in H}$, etc.
 
  Fix a function $\eps:H\to {\Bbb C}^*$ such that $\eps(\tau)+\eps(\bar{\tau})
  =0$ for any $\tau\in H$. 
  We define the symplectic form $\omega $ on $\xvw$ by
  $$ \omega ((B,t,s),(B',t',s'))=
         \hsum tr(\eps (\tau)B_{\bar{\tau}}B'_{\tau})+\isum 
  tr( s_i t'_i-s'_i t_i)\leqno(3.2.1)$$
    The algebraic group $G(\nu)=\iprod GL(V_i)$ acts on $\xvw$ by
  $$(B,t,s)\mapsto (g_{\inn(\tau)}B_{\tau}g_{\out(\tau)}^{-1},
  t_ig_i^{-1},g_is_i)\leqno(3.2.2)$$
  where $g=(g_i)\in G(\nu)$. The action of $G(\nu)$ preserve 
  the symplectic form
  $\omega$. Let $\mu :\xvw \to {\frak g}(\nu)$ be the associated moment map.
  Its $i$-th component $\mu_i:
  \xvw \to \End(V_i)$ is given by
  $$\mu_i((B,t,s))=\mathop{\sum}_{\tau\in H,i=\out(\tau)} 
  \eps(\tau)B_{\bar{\tau}}B_{\tau}+s_it_i.$$
 \subsection{}
 We recall basic properties of Nakajima's quiver variety. The results of this
 subsection originally proved by Nakajima [15],[16].
 \begin{defn}
  {\em We denoted by $\xsvw$ the set of all elements $(B,t,s)\in \xvw$ 
  satisfying the 
  following property; if there is a family of subspaces $V'_i$ in $V_i$ such 
  that they are  $B$-invariant ({\it i.e}, for each $\tau\in H$, 
  $B_{\tau}(V'_{\out(\tau)})\subset V'_{\inn(\tau)}$) 
  and contained to $\Ker(t_i)$ for each $i\in I$,then $V'_i=0$ for each 
  $i\in I$. We call an element of $\xsvw$ a stable point of $\xvw$}.
 \end{defn}
 It is easy to see that $\xsvw$ is an open subset of $\xvw$. 
 Clearly $G(\nu)$ acts on $\xsvw$.
 \begin{lemma}
 The action of $G(\nu)$ on $\xsvw$ is fixed point free.
 \end{lemma}
  By this lemma we can consider the quotient variety of $\xsvw$ by $G(\nu)$.
  That is, let
   $$\xxvw=(\mu^{-1}(0) \cap \xsvw) / G(\nu). $$
  We call $\xxvw$ {\it a quiver variety}. We denote by 
  $[B,t,s]$ the $G(\nu)$-orbit of $(B,t,s)$ considered as a point in $\xxvw$.
  \begin{prop} Assume $\xsvw\ne \phi$.
  \begin{enumerate}
   \arab
   \item $\xxvw$ is a quasi-projective smooth variety of dimension
         $\parallel \lambda\parallel^2-\\
         \parallel \lambda
         +\nu\parallel^2$. 
   \item  $\xxvw$ has the symplectic structure induced by $\omega$.
 \end{enumerate}
 \end{prop}
%
\section{Construction of crystal base}
 \subsection{}
 Let $\nu$, $\bar{\nu}\in Q_-$ such that $\nu-\bar{\nu}\in {\Bbb Z}_{\leq 0}
 \alpha_i$ for some $i\in I$. Assume that $V\in {\cal V}_{\nu}
 ,\bar{V}\in {\cal V}_{\bar{\nu}}$.
 Now we consider the diagram
  $$\xbvw\mathop{\longleftarrow}^{q_1}\xbvvw
  \mathop{\longrightarrow}^{q_2}\xvw.
 \leqno(4.1.1)$$
 Here $\xbvvw$ is the variety of $(B,t,s,\phi)$
 where $(B,t,s)\in \xvw$ and $\phi=(\phi_i):\bar{V}\to V$ is a injective
 morphism of $I$-graded vector spaces such that $\Imm\phi=(\Imm\phi_i)$ 
 is $B$-stable and contains $\Imm s=(\Imm s_i)$. Hence $B,t,s$ induce
 $\bar{B}:\bar{V}\to\bar{V}$, $\bar{t_i}:\bar{V_i}\to W_i$ and $\bar{s_i}:
 W_i\to\bar{V_i}$ respectively.
 We define $q_1(B,t,s,\phi)=(\bb,\bt,\bs)$ and
 $q_2(B,t,s,\phi)=(B,t,s)$.
 The following lemma is proved easily.
 \begin{lemma}
  We use the above notations.
  \begin{enumerate}
   \arab   
   \item The following two conditions are equivalent;
         \begin{enumerate}
          \renewcommand{\labelenumi}{(\alph{enumi})}
           \item $(B,t,s)$ is a stable point.
           \item $(\bb,\bt,\bs)$ is a stable point and the map 
                 $$V_i \xrightarrow{(B_{\tau},t_i)} 
                 \mathop{\oplus}_{\tau;\out(\tau)=i}V_{\inn(\tau)}\oplus W_i$$
                 is injective.
         \end{enumerate}
   \item The following two conditions are equivalent;\\
           (a) $\mu((B,t,s))=0$.
           (b) $\mu((\bb,\bt,\bs))=0$.
   \end{enumerate}  
  \end{lemma}
 By this lemma we can restrict the diagram to ${\mu}^{-1}(0)$ and stable
 points and take its quotient by $G(\nu),G(\bar{\nu})$. Then we get the diagram
 $$\xxbvw\mathop{\longleftarrow}^{\varpi_1}\xxbvvw
 \mathop{\longrightarrow}^{\varpi_2}\xxvw.\leqno(4.1.2)$$
 Here the variety of $([B,t,s],\Imm(\phi))$ where $[B,s,t]\in \xxvw$ and 
 $\phi$ is given above. 
\subsection{}
 For $i\in I$ and $c\in {\Bbb Z}_{\geq 0}$ we consider
 $$\xxvw_{i,c}=\{[B,t,s]\in \xxvw|\eps_i((B,t,s))=c\}$$
 where
 $$\eps_i((B,t,s))=\dimc \text{Coker}(\mathop{\oplus}_{\tau;\inn(\tau)=i}
 V_{\out(\tau)}\oplus W_i \xrightarrow{(B_{\tau},s_i)}V_i).$$
 Since $[B,t,s]$ is a $G(\nu)$-orbit the above definition is well defined.
 It is clear that $\xxvw_{i,c}$ is a locally closed 
 subvariety of $\xxvw$. 

 \begin{lemma}
 Suppose $\xxvw_{i,c}\ne \phi$. Then,\\
 $(1)$\quad $c+\langle h_i,\lambda+\nu\rangle\geq 0$.\\
 $(2)$\quad ${\frak X}(W;\nu-l\alpha_i)_{i,c+l}\ne \phi$ if and only if
 $-c\leq l \leq c+\langle h_i,\lambda+\nu\rangle$. 
 \end{lemma}
 {\it Proof.} First we shall show (1). Let $[B,t,s]\in \xxvw_{i,c}$. We 
 consider the following diagram;
 $$V_i \xrightarrow{(\varepsilon(\tau)B_{\tau},t_i)}  
   \mathop{\oplus}_{\tau;\out(\tau)=i}V_{\inn(\tau)}\oplus W_i
   \xrightarrow{(B_{\tau},s_i)}V_i.$$
 Since $(B,t,s)$ is a stable point the first map is injective. On the other
 hand $(B,t,s)\in \mu^{-1}(0)$ implies that the composition of these two 
 maps equal to zero. Therefore we have
 $$ \dimc V_i\leq \dimc \mbox{Coker}((B_{\tau},s_i))-\dimc V_i+
   \mathop{\sum}_{\tau;\out(\tau)=i}V_{\inn(\tau)}+\dimc W_i.$$
 Rewriting the above inequality, we get the statement.

 Let us prove the sufficient part of (2). Assume ${\frak X}
 (W;\nu-l\alpha_i)_{i,c+l}\ne \phi$. By (1) we have $c+l+\langle h_i,
 \lambda+\nu-l\alpha_i\rangle\geq 0$. Therefore we have
 $c+\langle h_i,\lambda+\nu\rangle\geq l$. By the definition of $\eps_i$
 we have $c+l\geq 0$. We get the statement.

  We shall show the necessary part of (2). Let $[B,t,s]\in \xxvw_{i,c}$.
 First we assume $-c\leq l\leq 0$ and let $\bar{\nu}=\nu-l\alpha_i$. We 
 introduce
 $\bar{V}\in {\cal V}_{\bar{\nu}}$ in the following way. If $j\ne i$, define 
 $\bar{V}_j
 =V_j$. $\bar{V}_i$ is $\langle h_i,\lambda+\bar{\nu}\rangle$-dimensional
 subspace of $V_i$ such that 
 $$\mbox{Im}(\mathop{\oplus}_{\tau;\inn(\tau)=i}V_{\out(\tau)}\oplus
 W_i\xrightarrow{(B_{\tau},s_i)}V_i)\subset \bar{V}_i.$$
 Let $(\bb,\bt,\bs)\in p_1p_2^{-1}(B,t,s)$. Then by Lemma 4.1.1 we have
 $(\bb,\bt,\bs)$ is a stable point and $(\bb,\bt,\bs)\in \mu^{-1}(0).$
 Moreover $\eps_i(\bb,\bt,\bs)=c+l$. 
 Therefore $[\bb,\bt,\bs]\in {\frak X}(W;\bar{\nu})_{i,c+l}$. That is
 ${\frak X}(W;\bar{\nu})_{i,c+l}\ne \phi$.

 Assume $0\leq l\leq c+\langle h_i,\lambda+\nu\rangle$. Let $\nu'=
 \nu-l\alpha_i$. We introduce $V'\in {\cal V}_{{\nu}'}$ in the 
 following way. If $j\ne i$ define 
 $V'_j
 =V_j$. $V'_i$ is $(l+\dimc V_i)$-dimensional subspace of $\mbox{Ker}
 (\mathop{\oplus}_{\tau;\inn(\tau)=i}V_{\out(\tau)}\oplus
 W_i\xrightarrow{(B_{\tau},s_i)}V_i).$ Note that the dimension of the kernel
 is $\langle h_i,\lambda+\nu\rangle+c+\dimc V_i$. Therefore such $V'$ exists. 
 Let $(B',t',s')\in p_2p_1^{-1}(B,t,s)$ such that the map
 $$V'_i\xrightarrow{(B'_{\tau},t'_i)}\mathop{\oplus}_{\tau;\out(\tau)=i}
   V_{\inn(\tau)}\oplus W_i\leqno{(4.2.1)}$$
 is injective. By the definition of $V'$ such $(B',t',s')$ exists.
 By Lemma 4.1.1 and (4.2.1) we have $(B',t',s')$ is a stable point and
 $\mu^{-1}(B',t',s')=0$. From the construction we have $\eps_i(B',t',s')=c+l$. 
 Therefore
 we conclude $[B',t',s']\in {\frak X}(W;\nu')_{i,c+l}$. That is ${\frak X}
 (W;\nu')_{i,c+l}\ne \phi$. 
 \hfill$\square$

 Assume $\xxvw_{i,c}\ne\phi$. By the above lemma we have $\xxbvw_{i,0}\ne\phi$
 for $\bar{\nu}=\nu+c\alpha_i$. Moreover, from the definition of the diagram 
 (4.1.2), we have 
 $$\varpi_1^{-1}(\xxbvw_{i,0})=\varpi_2^{-1}
  (\xxvw_{i,c}).$$
 We set
  $$\xxbvvw_{i,0}=\varpi_1^{-1}(\xxbvw_{i,0})=\varpi_2^{-1}(\xxvw_{i,c}).$$ 
 Then we have the following diagram
 $$\xxbvw_{i,0}\xleftarrow{\varpi_1}\xxbvvw_{i,0}\xrightarrow{\varpi_2}
  \xxvw_{i,c}\leqno(4.2.2).$$
 It is easy to see that the restriction of $\varpi_2$ to
 $\xxbvvw_{i,0}$ is an isomorphism and $\xxbvw_{i,0}$ is a open subvariety 
 of $\xxbvw$.
 \begin{lemma}
  \begin{enumerate}
  \arab
   \item For any $i\in I$,
         $${\frak X}(W;0)_{i,c}=
         \begin{cases}
         \mbox{\{pt.\}},& c=0,\\
         \phi,& c>0.
         \end{cases}$$
   \item Suppose $\xxvw_{i,c}\ne \phi$ and let $\bar{\nu}=\nu+c\alpha_i$. 
         Consider the restriction of 
         $\varpi_1$ to $\xxbvvw_{i,0}$. Then the
         fiber of this map is isomorphic to $\Grass_c({\Bbb C}^{\langle h_i,
         \lambda+\bar{\nu}\rangle})$.
         Here we denote by $\Grass_c({\Bbb C}^p)$ the Grassmanian variety of 
         $c$-dimensional subspaces of ${\Bbb C}^p$.  
  \end{enumerate}
 \end{lemma}
 {\it Proof.}
 We have (1) immediately from the definition. Let 
 $[\bar{B},\bar{t},\bar{s}]\in \xxbvw_{i,0}$ and $[B,t,s,\phi]\in 
 \xxbvvw_{i,0}$ such that $q_1((B,t,s,\phi))=(\bar{B},\bar{t},\bar{s})$. By
 the definition we have the following commutative diagram;     
 $$\bar{V}_i \xrightarrow{(\varepsilon(\tau)\bar{B}_{\tau},\bar{t}_i)}  
   \left(\mathop{\oplus}_{\tau;\out(\tau)=i}\bar{V}_{\inn(\tau)}\right)
   \oplus W_i\xrightarrow{(\bar{B}_{\tau},\bar{s}_i)}\bar{V}_i$$
 $$\hspace*{3mm}
   \downarrow{\mbox{\scriptsize{$\phi$}}}\hspace{35mm}
   \downarrow{\mbox{\scriptsize{$\phi\oplus \mbox{id}_{W_i}$}}}\hspace{20mm}
   \downarrow{\mbox{\scriptsize{$\phi$}}}$$
 $$V_i \xrightarrow{(\varepsilon(\tau)B_{\tau},t_i)}  
   \left(\mathop{\oplus}_{\tau;\out(\tau)=i}V_{\inn(\tau)}\right)\oplus W_i
   \xrightarrow{(B_{\tau},s_i)}V_i.$$
 We remark that the second vertical map is an isomorphism. By the similar 
 argument in the proof of lemma 4.2.1 we have the map 
 $(\varepsilon(\tau)B_{\tau},t_i)$ is injective and $\sum\eps(\tau)
 B_{\bar{\tau}}B_{\tau}+s_it_i=0$. On the other hand the map $(\bar{B}_{\tau},
 \bar{s}_i)$ is subjective because $[\bar{B},\bar{t},\bar{s}]$ is the element
 of $\xxbvw_{i,0}$. Therefore, for given $[\bar{B},\bar{t},\bar{s}]$, the
 fiber of $\varpi_1$ is the set of all elements $[B,t,s,\phi]\in 
 \xxbvvw_{i,0}$ such that
 $$\bar{V}_i\overset{\phi_i}{\hookrightarrow}V_i\hookrightarrow
   \mbox{Ker}\left(\mathop{\oplus}_{\tau;\out(\tau)=i}V_{\inn(\tau)}
   \twoheadrightarrow \bar{V}_i\right).$$
 Hence we have $\varpi_1^{-1}([\bar{B},\bar{t},\bar{s}])$ is isomorphic to
 $$\Grass_c({\Bbb C}^{\dimc W_i-\sum_{i\ne j}\langle \alpha_i,\alpha_j\rangle
   \dimc\bar{V}_j-2\dimc\bar{V}_i})=
 \Grass_c({\Bbb C}^{\langle h_i,
         \lambda+\bar{\nu}\rangle}).$$
 \hspace*{120mm}\hfill$\square$
 
 By Lemma 4.2.2 (2) we have the following corollary immediately.
 \begin{cor}Suppose $\xxvw_{i,c}\ne \phi$.
   There is one to one correspondence between the set of all 
          irreducible components of $\xxbvw_{i,0}$ and the set of all 
          irreducible components of $\xxvw_{i,c}$.
 \end{cor}
\subsection{}
 An element of $B=(B_{\tau})$ is said to be 
 nilpotent if there exists an
 integer $n\geq 2$ such that the following holds; for any sequence 
 $\tau_1,\tau_2,\cdots,\tau_n$ in $H$ such that $\inn(\tau_1)=\out(\tau_2),
 \inn(\tau_2)=\out(\tau_3),\cdots,\inn(\tau_{n-1})=\out(\tau_n)$, the 
 composition $B_{\tau_n}\cdots B_{\tau_2}B_{\tau_1}:V_{\out(\tau_1)}\to
 V_{\inn(\tau_n)}$ is zero.
 Let us define a subvariety $\lvw$ of $\xxvw$ by
 $$\lvw=\{[B,t,s]\in \xxvw|s=0\text{ and $B$ is nilpotent}\}.$$

 The following is due to Nakajima [15].
 \begin{prop}
 $\lvw$ is a Lagrangian subvariety of $\xxvw$.
 \end{prop} 
 Now we denote by $B(W;\nu)$ the set of all irreducible components of 
 $\Lambda(W;\nu)$. Take $\Lambda \in B(W;\nu)$. For a generic point $[B,t,s]$
 of $\Lambda$ we define $\eps_i(\Lambda)=\eps_i((B,t,s))$. For $c\in 
 {\Bbb Z}_{\geq 0}$ we set $B(W;\nu)_{i,c}$
 the set of all elements of $B(W;\nu)$ such that $\eps_i(\Lambda)=c$.
 \begin{prop}
 $$B(W;\bar{\nu})_{i,0}\cong B(W;\nu)_{i,c}.$$
 \end{prop}
{\it Proof.}
 We use the notations of Lemma 4.1.1. It is easy to see that $s=0$ if and only
 if $\bar{s}=0$ and $B$ is nilpotent if and only if $\bb$ is nilpotent. 
 By these equivalences we can restrict diagram (4.1.2) to nilpotent elements 
 and $s=0$. Therefore we obtain the desired result from Corollary 4.2.3.
\hfill$\square$
\subsection{}
 \begin{defn}
 {\em  Suppose that $\bar{\Lambda}\in B(W;\bar{\nu})_{i,0}$ corresponds to 
  $\Lambda\in 
  B(W;\nu)_{i,c}$ by the isomorphism in Proposition 4.3.2. Then we define maps 
  $\fit^c:B(W;\bar{\nu})_{i,0}
  \to B(W;\nu)_{i,c}$ and $\eit^c:B(W;\nu)_{i,c}\to B(W;\bar{\nu})_{i,0}$ by
  $$\fit^c(\bar{\Lambda})=\Lambda,$$
  $$\eit^c(\Lambda)=\bar{\Lambda}.$$
  Furthermore we define the maps 
  $$\eit,\text{} \fit:\mathop{\bigsqcup}_{\nu}B(W;\nu)\to
  \mathop{\bigsqcup}_{\nu}B(W;\nu)\sqcup\{0\}$$
  by
  $$\eit:B(W;\nu)_{i,c}\xrightarrow{{\eit}^c} B(W;\nu+c\alpha_i)_{i,0}
         \xrightarrow{{\fit}^{c-1}} B(W;\nu+\alpha_i)_{i,c-1},$$
  $$\fit:B(W;\nu)_{i,c} \xrightarrow{{\eit}^c} B(W;\nu+c\alpha_i)_{i,0}
         \xrightarrow{{\fit}^{c+1}} B(W;\nu-\alpha_i)_{i,c+1}$$
  We put $\eit(\Lambda)=0$ for $\Lambda\in B(W;\nu)_{i,0}$ and
  $\fit(\Lambda)=0$ for $\Lambda\in B(W;\nu)_{i,c}$ for $c<-\langle h_i,
  \lambda+\nu\rangle$}.
 \end{defn}
 Then the maps $\eit^c$ (resp. $\fit^c$) which is constructed in 
 the definition may be considered as a $c$-th power of $\eit$ (resp. $\fit$).
 Let us define a map $wt:
 \mathop{\bigsqcup}_{\nu}B(W;\nu) \to P$ by $wt(\Lambda)=\lambda+\nu \in P$ 
 for $\Lambda \in B(W;\nu)$. We set $\varphi_i(\Lambda)=\eps(\Lambda)+
 \langle h_i,wt(\Lambda)\rangle$.
 \begin{thm}
 $\mathop{\bigsqcup}_{\nu}B(W;\nu)$ is a crystal in the sense of Definition
 2.1.1.
 \end{thm}
 {\it Proof.} It is enough to see that the axioms of the crystal are
 satisfied. But it is obvious from the definition.\hfill$\square$
 \begin{lemma}
  For any $\Lambda\in \mathop{\bigsqcup}_{\nu}B(W;\nu)$, $\varphi_i(\Lambda)
  =\mbox{max}\{k\geq 0|\fit^k(\Lambda)\ne 0\}$.
 \end{lemma}
 {\it Proof.} Suppose $B(W;\nu)\ne \phi$ and $\Lambda\in B(W;\nu)_{i,c}$.
 By the definition $\varphi_i(\Lambda)=c+\langle h_i,\lambda+\nu\rangle$.
 By Lemma 4.2.1 $B(W;\nu-l\alpha_i)_{i,c+l}\ne \phi$ if and only if
 $-c\leq l\leq c+\langle h_i,\lambda+\nu\rangle$. Therefore we have the
 lemma.   
 \hfill$\square$
\subsection{}
 We recall the result of [7]. See [7] for details.

 Let $V\in {\cal V}_{\nu}$.
 Define
 $$ X(\nu) =(\hopl \Homc (\vout ,\vin))$$
 and the symplectic form $\omega$ on $X(\nu)$ by 
 $$ \omega (B,B')=
         \hsum tr(\eps (\tau)B_{\bar{\tau}}B).$$
 The algebraic group $G(\nu)$ acts on $X(\nu)$ in the same way to (3.2.2). 
 Let $\mu$ be the corresponding moment map. For $B\in X(\nu)$ we set
 $$\eps_i(B)=\dimc \text{Coker}(\mathop{\oplus}_{\tau;\inn(\tau)=i}
 V_{\out(\tau)}\xrightarrow{(B_{\tau})}V_i).$$
 Let 
 $$\Lambda(\nu)=\{B\in X(\nu)|\mu(B)=0\mbox{ and $B$ is nilpotent}\}$$
 and $B(\infty;\nu)$ the set of irreducible components of $\Lambda(\nu)$.
 For a generic point $B$ of $\Lambda\in B(\infty;\nu)$ we set $\eps_i(\Lambda)
 =\eps_i(B)$. We can define the operators
 $$\eit:\mathop{\bigsqcup}_{\nu}B(\infty;\nu)\to
  \mathop{\bigsqcup}_{\nu}B(\infty;\nu)\sqcup\{0\}, \quad
  \fit:\mathop{\bigsqcup}_{\nu}B(\infty;\nu)\to
  \mathop{\bigsqcup}_{\nu}B(\infty;\nu)$$
 by the similar way to Definition 4.4.1. Set $\mbox{wt}(\Lambda)=\nu$ and
 $\varphi_i(\Lambda)=\eps_i(\Lambda)+\langle h_i,\mbox{wt}(\Lambda)\rangle$
 for $\Lambda\in B(\infty;\nu)$.
 \begin{thm}
  \begin{enumerate}
   \arab
   \item $\mathop{\bigsqcup}_{\nu}B(\infty;\nu)$ is a crystal.
   \item $\mathop{\bigsqcup}_{\nu}B(\infty;\nu)$ is isomorphic to $B(\infty)$.
  \end{enumerate}
 \end{thm}
\subsection{ }
 From Proposition 4.3.1 we have the following statement.
 \begin{lemma}
  $B(W;\nu)$ is isomorphic to the set of irreducible 
  $G(\nu)$-invariant Lagrangian subvarieties of $\xsvw$.
 \end{lemma}
 Define 
 $$\tilde{\Lambda}(W;\nu)=\{(B,t,s)\in \xvw|s=0\mbox{ and $B$ is nilpotent}\}$$
 and let $\tilde{B}(W;\nu)$ be the set of all irreducible components of 
 $\tilde{\Lambda}(W;\nu)$. From Lemma 4.6.1 there exist an injective map
 $$\iota: B(W;\nu)\to \tilde{B}(W;\nu).$$
 We introduce the map $\kappa: \tilde{\Lambda}(W;\nu)\to \Lambda(\nu)$ by
 $(B,t,s)\mapsto B$. It is easy to see that $\kappa$ induce the following
 isomorphism;
 $$\kappa:\tilde{B}(W;\nu)\cong B(\infty;\nu).$$ 
 Therefore we have the following.
 \begin{lemma}
  There is an injective map $\kappa\circ\iota:B(W;\nu)\to B(\infty;\nu)$.
 \end{lemma}  
 Let us introduce the map $\Psi:T_{\lambda}\otimes\mathop{\bigsqcup}_{\nu}
 B(\infty;\nu)\to 
 \mathop{\bigsqcup}_{\nu}B(W;\nu)\sqcup\{0\}$ by
 $$
 \Psi(t_{\lambda}\otimes\Lambda)=
  \begin{cases}
   (\kappa\circ\iota)^{-1}(\Lambda),\quad & \mbox{if } \Lambda\in\mbox{Im}
   (\kappa\circ\iota),\\
   0, &\mbox{otherwise.}
  \end{cases}
 $$
 \begin{lemma}
  $\Psi$ is a strict morphism of crystal with the following properties.
  \begin{enumerate}
   \arab
   \item $\mbox{Im}\Phi=\mathop{\bigsqcup}_{\nu}B(W;\nu)\sqcup\{0\}$.
   \item The set $\{t_{\lambda}\otimes \Lambda\in (T_{\lambda}\otimes
         \mathop{\bigsqcup}_{\nu}B(\infty;\nu))~|~\Psi(t_{\lambda}\otimes
         \Lambda)\ne 0\}$ is isomorphic to $B(\lambda)$ through $\Phi$.
  \end{enumerate}
 \end{lemma}
 {\it Proof.} By the definition of $\Phi$ (1) and (2) are clear. Therefore
 it is enough to show that $\Phi$ is a strict morphism. From the construction
 $\Phi$ preserves weights. Let $\Lambda\in B(W;\nu)$ and $[B,t,s]\in \Lambda$.
 Since $s=0$ and $\eps_i(t_{\lambda})=\varphi_i(t_{\lambda})=-\infty$ the map 
 $\kappa\circ\iota$ preserves the values of $\eps_i$
 and $\varphi_i$. Therefore $\Phi$ also preserves them. By the definition
 $\kappa\circ\iota$ commutes with $\eit$ and $\fit$. This means $\Phi$
 commutes with them.   
 \hfill$\square$
 \begin{thm}
  $\mathop{\bigsqcup}_{\nu}B(W;\nu)$ is isomorphic to $B(\lambda)$ as a 
  crystal.
 \end{thm}
 {\it Proof.} By Theorem 4.5.1 we can consider $\Phi$ is a strict morphism
 from $T_{\lambda}\otimes B(\infty)$ to $\mathop{\bigsqcup}_{\nu}B(W;\nu)$.
 It is enough to see that the 
 conditions in Proposition 2.3.1 are satisfied. (1) follows from
 the definition. (2) and (4) are just the above lemma. (3) is already
 proved in Lemma 4.4.3. \hfill$\square$
\subsection{ }
 Recall that $\lvw$ is a Lagrangian subvariety of $\xxvw$. It is homotopic 
 to $\xxvw$ (See [16].), hence the top 
 homology group of $\lvw$ is isomorphic to the middle 
 homology group of $\xxvw$. Since the set of all irreducible components of
 $\lvw$ is isomorphic to $B(\lambda)_{\lambda+\nu}$, we have
 $$\dim H_{\mbox{\tiny{middle}}}(\xxvw)=\sharp B(\lambda)_{\lambda+\nu}.$$
 Here $\nu=-\sum\dim V_i\alpha_i\in Q_-$ and $B(\lambda)_{\lambda+\nu}$ 
 is the set of all elements of $B(\lambda)$ with weight $\lambda+\nu$.
 As an application we can compute the dimension of the middle homology
 groups of quiver varieties by using Kac's character formula [2]. That is, 
 $$\mathop{\sum_{\nu\in Q_-}}\dim H_{\mbox{\tiny{middle}}}(\xxvw)
 e^{\lambda+\nu}=
 \frac{\mathop{\sum_{w\in W}}(-1)^{l(w)}e^{w(\lambda+\rho)-\rho}}
 {\mathop{\prod_{\alpha\in\Delta_+}}(1-e^{-\alpha})^{\mbox{\tiny{mult}}
 \alpha}}.$$
 Here $W$ is the Weyl group of $\frak g$, $l(w)$ the length of $w\in W$,
 $\rho$ a Weyl vector, $\Delta_+$ the set of all positive roots and 
 $\mbox{mult}\alpha$ the multiplicity of a root $\alpha$.  
 This result was already given by Lusztig and Nakajima. (See [12]
 and [16].) We can reprove it in a crystal-theoretical way.
 
\mbox{ }\\
 \begin{center}
  {\small Department of Mathematics, Hiroshima University,\\
  Higashi-Hiroshima 739-8526,Japan}\\
  {\small e-mail : {\tt yosihisa@math.sci.hiroshima-u.ac.jp}}
 \end{center}
\end{document}